\documentclass{conm-p-l}
\usepackage[all]{xy}
\usepackage{amssymb}%,MnSymbol}
\usepackage[hypertex]{hyperref}
\usepackage{graphicx}
\newtheorem{theorem}{Theorem}[section]
\newtheorem{lemma}[theorem]{Lemma}
\newtheorem{prop}[theorem]{Proposition}
\theoremstyle{definition}
\newtheorem{definition}[theorem]{Definition}
\newtheorem{example}[theorem]{Example}
\newtheorem{xca}[theorem]{Exercise}

\theoremstyle{remark}
\newtheorem{remark}[theorem]{Remark}

\numberwithin{equation}{section}
\newcommand{\C}{\mathbb C}
\newcommand{\R}{\mathbb R}  %@@
\newcommand{\indhg}{\operatorname{ind}_{\mathfrak h}^{\mathfrak g}}
\newcommand{\indpg}{\operatorname{ind}_{\mathfrak p}^{\mathfrak g}}
\newcommand{\indpgprime}{\operatorname{ind}_{\mathfrak p'}^{\mathfrak g'}}
\newcommand{\LieG}{\mathrm{Lie}(G)}
\newcommand{\To}{\longrightarrow}
\newcommand{\abs}[1]{\lvert#1\rvert}

\newcommand{\blankbox}[2]{%
  \parbox{\columnwidth}{\centering
    \setlength{\fboxsep}{0pt}%
    \fbox{\raisebox{0pt}[#2]{\hspace{#1}}}%
  }%
}

\begin{document}

\title{$F$-method for constructing equivariant differential operators%
}

\author{Toshiyuki KOBAYASHI}
\address{Kavli IPMU (WPI), and Graduate School of Mathematical Sciences \\
The University of Tokyo}
%\curraddr{Department of Mathematics and Statistics,
%Case Western Reserve University, Cleveland, Ohio 43403}
\email{toshi@ms.u-tokyo.ac.jp}
\thanks{Partially supported by IHES and Grant-in-Aid for Scientific Research 
(B) (22340026) JSPS}

%\author{Author Two}
%\address{Mathematical Research Section, School of Mathematical Sciences,
%Australian National University, Canberra ACT 2601, Australia}

%    General info
\date{
October 16, 2012.
}

\dedicatory{Dedicated to  Professor Sigurdur Helgason
 for his 85th birthday.}


\begin{abstract}
Using an algebraic Fourier transform of operators,
we develop a method ($F$-method) to obtain explicit highest weight vectors
in the branching laws
 by differential equations.
This article gives a brief explanation of the $F$-method
and its applications 
 to a concrete construction of some natural equivariant operators that arise in
parabolic geometry and in automorphic forms.
\end{abstract}

\maketitle

\setcounter{tocdepth}{1}
\tableofcontents

\section{Introduction}
The aim of this article is to give a brief account of a method
 that helps us to find a closed formula of highest weight vectors
in the branching laws of certain generalized Verma modules,
or equivalently,
to construct explicitly equivariant differential operators
from generalized flag varieties to subvarieties.

This method, which we call the $F$-method,
 transfers an algebraic problem
 of finding explicit highest weight vectors
 to an analytic problem
 of solving differential equations
 (of second order)
 via the algebraic Fourier transform of operators
(Definition \ref{def:hat}). 
A part of the ideas of the $F$-method has grown in a detailed analysis of
the Schr\"{o}dinger model of the minimal
representation of indefinite orthogonal groups
\cite{intopq}.

The $F$-method provides a conceptual understanding of some
 natural differential
operators 
which were previously found by a combinatorial approach based on 
recurrence formulas.
Typical examples that we have in mind are
the Rankin--Cohen bidifferential operators 
\[
R_n^{k_1,k_2} (f_1,f_2)
= \sum_{j=0}^n (-1)^j \begin{pmatrix} n \\ j \end{pmatrix}
   \frac{(k_1+n-1)!(k_2+n-1)!}{(k_1+n-j-1)!(k_2+j-1)!}
   \frac{\partial^{n-j}f_1}{\partial x^{n-j}}
   \frac{\partial^j f_2}{\partial y^j}
   \bigg|_{x=y}
\]
in  
automorphic form theory 
\cite{CMZ,vDP,Ra},
and Juhl's conformally equivariant operators \cite{Juhl}
from $C^\infty(\mathbb{R}^n)$ to $C^\infty(\mathbb{R}^{n-1})$:
\[
T_{\lambda,\nu}
= \sum_{2j+k=\nu-\lambda}
   \frac{1}
          {2^j j!(\nu-\lambda-2j)!}
   \biggl(\prod_{i=1}^{\frac{\nu-\lambda}{2}-j} (\lambda+\nu-n-1+2i) \biggr)
   \Delta_{\mathbb{R}^{n-1}}^j \Bigl( \frac{\partial}{\partial x_n} \Bigr)^k .
\]

These examples can be reconstructed by the $F$-method by using 
a special case
 of the {\textit{fundamental differential operators}}, 
 which are commuting family
 of second order differential operators
 on the isotropic cone,
 see \cite[(1.1.3)]{intopq}.

In recent joint works
with
B.~{\O}rsted, M.~Pevzner, P.~Somberg and V.~Sou\v{c}ek
\cite{KOSS,KP},
we have developed the $F$-method
 to more general settings,
 and have found new explicit formulas of
equivariant differential operators in parabolic geometry,
and also have obtained a generalization of the Rankin--Cohen operators. 
To find those nice settings
 where the $F$-method works well,
 we can apply the general theory 
  \cite{K08,K10}
that assures discretely decomposable
 and multiplicity-free restrictions
 of representations
 to reductive subalgebras.  

The author expresses his sincere gratitude to the organizers,
J. Christensen, F. Gonzalez and T. Quinto,
for their warm hospitality
 during the conference
 in honor of Professor Helgason
 for his 85th birthday in Boston 2012.
The final manuscript was prepared when the author was visiting IHES.

\section{Preliminaries}
\subsection{Induced modules}
Let $\mathfrak g$ be a Lie algebra over $\C$, and $U(\mathfrak g)$ its universal enveloping algebra.
Suppose that $\mathfrak h$ is a subalgebra of $\mathfrak g$
 and $V$ is an $\mathfrak h$-module.
We define the induced $U(\mathfrak g)$-module by
$$
\indhg(V):=U(\mathfrak g)\otimes_{U(\mathfrak h)}V.
$$
If $\mathfrak h$ is a Borel subalgebra and if
$\dim_{\mathbb{C}}V = 1$, then $\indhg(V)$ is the standard Verma module.

\subsection{Extended notion of differential operators}

We understand clearly the notion of differential operators
between two vector bundles over the same base manifold.
We extend this notion in a more general setting
where there is a morphism between two base manifolds
(see \cite{KP} for details).
\begin{definition}\label{def:21}
Let $\mathcal V_X\to X$ and $\mathcal W_Y\to Y$ be two vector bundles 
with a smooth map  $p:Y\to X$ between the base manifolds.
Denote by $C^\infty(X,\mathcal{V}_X)$ and
$C^\infty(Y,\mathcal{W}_Y)$
the spaces of smooth sections to the vector bundles.
We say that a linear map $T:C^\infty(X,\mathcal V_X)\to C^\infty(Y,\mathcal W_Y)$ is a \emph{differential
operator} 
if $T$ is a local operator in the sense that 
\begin{equation}
\label{eqn:suppT}
\operatorname{Supp}(T f)
 \subset
 {p^{-1}(\operatorname{Supp}f)},\qquad \mathrm{for\, any}\,
 f\in C^\infty (X,\mathcal V).
\end{equation}
We write $\operatorname{Diff}(\mathcal V_X,\mathcal W_Y)$ for the vector space of such differential operators.
\end{definition}

Since any smooth map $p:Y \to X$ is given
 as the composition of a submersion and an embedding
\[
   Y \hookrightarrow X \times Y \twoheadrightarrow X,
  \quad
  y \mapsto (p(y), y) \mapsto p(y), 
\] 
the following example describes the general
situation.
\begin{example}
Let $n$ be the dimension of $X$.  
\begin{enumerate}
\item 
Suppose $p:Y\twoheadrightarrow X$ is a submersion. 
Choose local coordinates $\{(x_i,z_j)\}$ on $Y$ 
 such that $X$ is given locally
 by $z_j=0$. 
Then every $T\in \operatorname{Diff}(\mathcal V_X,\mathcal W_Y)$
 is locally of the form
$$
\sum_{\alpha\in{\mathbb{N}}^n}h_\alpha(x,z)\frac{\partial^{|\alpha|}}{\partial x^\alpha},
$$
where $h_\alpha(x,z)$ are $\operatorname{Hom}(V,W)$-valued smooth functions on $Y$.

\item
Suppose $i: Y\hookrightarrow X$ is an embedding. 
Choose local coordinates $\{(y_i,z_j)\}$ on $X$ 
 such that $Y$ is given locally by $z_j=0$. 
Then every $T\in \operatorname{Diff}(\mathcal V_X,\mathcal W_Y)$
 is locally of the form
$$
\sum_{(\alpha,\beta)\in{\mathbb{N}}^n}g_{\alpha\beta}(y)\frac{\partial^{|\alpha|+|\beta|}}{\partial y^\alpha\partial z^\beta},
$$
where $g_{\alpha,\beta}(y)$ are $\operatorname{Hom}(V,W)$-valued smooth functions on $Y$.
\end{enumerate}
\end{example}

\subsection{Equivariant differential operators}\label{sec:33}

Let $G$ be a real Lie group, $\mathfrak g(\R)=\LieG$ and $\mathfrak g=
\mathfrak g(\R)\otimes\C$. Analogous notations will be applied to other
Lie groups denoted by uppercase Roman  letters.

Let $dR$ be the representation of $U(\mathfrak g)$ on the space 
$C^\infty (G)$ of smooth complex-valued functions on $G$
 generated by the Lie algebra action:
\begin{equation}\label{eqn:dR}
 (dR(A)f)(x):=\frac{d}{dt}\Big\vert_{t=0}f(xe^{tA})
\quad\text{for 
$A \in \mathfrak{g}(\mathbb{R})$}.
\end{equation}

Let $H$ be a closed subgroup of $G$. Given a finite dimensional representation $V$  of $H$ we form a homogeneous  vector bundle 
$\mathcal V_X:= G\times_HV$ over the homogeneous space $X:=G/H$.
The space of smooth sections $C^\infty(X,\mathcal V_X)$ can be seen
as a subspace of $C^\infty(G)\otimes V$.

Let $V^\vee$ be the (complex linear) dual space of $V$.
Then the $(G \times \mathfrak{g})$-invariant bilinear map
$
C^\infty(G)\times U(\mathfrak g)\to C^\infty(G)$,
$
(f,u)\mapsto dR(u)f$
induces a commutative diagram of $(G \times \mathfrak{g})$-bilinear maps:
$$
\begin{matrix}
C^\infty(G)\otimes V\times U(\mathfrak g)\otimes_\C V^\vee
&\longrightarrow
&C^\infty(G)
\\
\text{\rotatebox{90}{$\hookrightarrow$}} 
\qquad
\text{\rotatebox{90}{$\twoheadleftarrow$}} 
&
&
% \lhookuparrow\qquad\twoheaddownarrow
\parallel
\\
C^\infty(X,\mathcal V_X)\times \indhg(V^\vee)
&
\longrightarrow&
C^\infty(G).
\end{matrix}
$$

In turn, we get the following natural $\mathfrak{g}$-homomorphism:
\begin{equation}\label{eqn:1234}
\indhg(V^\vee)\To\operatorname{Hom}_G(C^\infty(X,\mathcal V_X),
C^\infty(G)).
\end{equation}

Next, we take a connected closed subgroup $H'$ of $H$.
For a
 finite dimensional representation $W$ of $H'$ we form the
homogeneous vector bundle $\mathcal W_Z:=G\times_{H'}W$ over $Z:=G/H'$. 
Taking the tensor product of \eqref{eqn:1234} with $W$,
and collecting all $\mathfrak{h}'$-invariant elements,
we get an injective homomorphism:
\begin{equation}\label{eqn:morf}
\operatorname{Hom}_{\mathfrak h'}(W^\vee,\indhg(V^\vee))\To
\operatorname{Hom}_G(C^\infty(X,\mathcal V_X),
C^\infty(Z,\mathcal W_Z)),\quad \varphi\mapsto D_\varphi.
\end{equation}
 
Finally,
 we take any closed subgroup $G'$ containing $H'$ and form a homogeneous vector bundle $\mathcal W_Y:= G'\times_{H'} W$
 over $Y:=G'/H'$.
We note that $\mathcal{W}_Y$ is obtained from $\mathcal{W}_Z$
by restricting the base manifold $Z$ to $Y$. 

Let $R_{Z\to Y}:C^\infty(Z,\mathcal W_Z)\to C^\infty(Y,\mathcal W_Y)$ be the restriction map.
We set
\begin{equation}\label{eqn:DXY}
D_{X\to Y}(\varphi):=R_{X\to Y}\circ D_\varphi.
\end{equation}
Since there is a natural ($G'$-equivariant but not necessarily injective) morphism
$Y \to X$,
the extended notion of differential operators between
$\mathcal{V}_X$ and $\mathcal{W}_Y$ makes sense
(see Definition \ref{def:21}).
We then have:
\begin{theorem}\label{thm:surject}
The operator $D_{X\to Y}$ (see \eqref{eqn:DXY}) induces a bijection:
\begin{eqnarray}\label{eqn:isoo}
D_{X\to Y}:\operatorname{Hom}_{\mathfrak h'}(W^\vee,\indhg(V^\vee))
\stackrel{\sim}{\longrightarrow}
\operatorname{Diff}_{G'}( \mathcal V_X,\mathcal W_Y).
\end{eqnarray}
\end{theorem}

\begin{remark}\label{rem:cpx}
We may consider a holomorphic version of Theorem~\ref{thm:surject}
 as follows.
Suppose $G_{\mathbb{C}}$, $H_{\mathbb{C}}$, $G'_{\mathbb{C}}$ 
and $H'_{\mathbb{C}}$ are connected complex Lie groups
with Lie algebras $\mathfrak{g}$, $\mathfrak{h}$, $\mathfrak{g}'$ and
$\mathfrak{h}'$, 
and $\mathcal{V}_{X_{\mathbb{C}}}$ and
$\mathcal{W}_{Y_{\mathbb{C}}}$ are homogeneous
holomorphic vector bundles over
$X_{\mathbb{C}} := G_{\mathbb{C}} / H_{\mathbb{C}}$
and
$Y_{\mathbb{C}} := G'_{\mathbb{C}} / H'_{\mathbb{C}}$,
respectively.
Then Theorem \ref{thm:surject}
 implies that we have a bijection:
\begin{equation}\label{eqn:CDXY}
D_{X_{\mathbb{C}} \to Y_{\mathbb{C}}} :
\operatorname{Hom}_{\mathfrak{h}'}
(W^\vee, \operatorname{ind}_{\mathfrak{h}}^{\mathfrak{g}} (V^\vee))
\overset{\sim}{\longrightarrow}
\operatorname{Diff}_{G'_{\mathbb{C}}}^{\operatorname{hol}}
(\mathcal{V}_{X_{\mathbb{C}}}, \mathcal{W}_{Y_{\mathbb{C}}}).
\end{equation}
Here
$\operatorname{Diff}_{G'_{\mathbb{C}}}^{\operatorname{hol}}$
denotes the space of $G'_{\mathbb{C}}$-equivariant holomorphic differential operators with respect to the holomorphic map $Y_{\mathbb C} \to X_{\mathbb C}$.
By the universality of the induced module,
\eqref{eqn:CDXY} may be written as
\begin{equation}\label{eqn:CDXY2}
D_{X_{\mathbb{C}} \to Y_{\mathbb{C}}} :
\operatorname{Hom}_{\mathfrak{g}'}(\operatorname{ind}_{\mathfrak{h}'}^{\mathfrak{g}'}(W^\vee),
\operatorname{ind}_{\mathfrak{h}}^{\mathfrak{g}} (V^\vee))
\overset{\sim}{\longrightarrow}
\operatorname{Diff}_{G'_{\mathbb{C}}}^{\operatorname{hol}}
(\mathcal{V}_{X_{\mathbb{C}}}, \mathcal{W}_{Y_{\mathbb{C}}}).
\end{equation}
\end{remark}

The isomorphism \eqref{eqn:CDXY2} is well-known
when $X_{\mathbb{C}}=Y_{\mathbb{C}}$
 is a complex flag variety.  
The
  proof of Theorem \ref{thm:surject} 
is given in  \cite{KP}
in the generality
that $X \ne Y$.

\subsection{Multiplicity-free branching laws}
Theorem \ref{thm:surject} says that if 
$\operatorname{Hom}_{\mathfrak{h}'} (W^\vee,\operatorname{ind}_{\mathfrak{h}}^{\mathfrak{g}}(V^\vee))$
is one-dimensional
 then 
 $G'$-equivariant differential operators from
$\mathcal{V}_X$ to $\mathcal{W}_Y$ 
are unique up to scalar.
Thus we may expect that such unique operators
should have a natural meaning
and would be given by a reasonably simple formula.
Then we may be interested in finding systematically 
the examples where 
$\operatorname{Hom}_{\mathfrak{h}'} (W^\vee,\operatorname{ind}_{\mathfrak{h}}^{\mathfrak{g}}(V^\vee))$
is one-dimensional.
This is a special case of the \textit{branching problems} 
that asks how representations
decompose when restricted to subalgebras. 
In the setting where  $\mathfrak{h}$ is a parabolic subalgebra 
(to be denoted by
${\mathfrak{p}}$)
of  a reductive Lie algebra $\mathfrak{g}$,
we have the following theorem:
\begin{theorem}\label{thm:1.4}
Assume the nilradical\/ $\mathfrak{n}_+$ of\/ $\mathfrak{p}$ is abelian and $\tau$ is an involutive
automorphism of\/ $\mathfrak{g}$ such that
$\tau\mathfrak{p} = \mathfrak{p}$.
Then for any one-dimensional representation
$\mathbb{C}_\lambda$ of\/ $\mathfrak{p}$ 
and for any finite dimensional representation $W$ of\/
$\mathfrak{p}^\tau :=\{ X \in \mathfrak{p}: \tau X = X \}$,
we have
\[
\dim \operatorname{Hom}_{\mathfrak{p}^\tau}
	   (W^\vee, \operatorname{ind}_{\mathfrak{p}}^{\mathfrak{g}}
    (\mathbb{C}_\lambda^\vee)) \le 1.
\]
\end{theorem}
There are two known approaches for the proof of Theorem \ref{thm:1.4}.
One is geometric --- to use the general theory of
 the \textit{visible action} on complex manifolds \cite{visible,K08},
and the other is algebraic --- to work inside the universal enveloping algebra \cite{K10}.
\begin{remark}
\label{rem:Branching}
Branching laws in the setting of Theorem \ref{thm:1.4}
 are explicitly obtained
 in terms of \lq{relative strongly orthogonal roots}\rq\
on the level of the Grothendieck group,
which becomes a direct sum decomposition when the parameter
$\lambda$ of $V$ is `generic' or sufficiently positive, 
\cite[Theorems 8.3 and 8.4]{K08} or \cite{K10}.  
The $F$-method will give a finer structure
 of branching laws 
 by finding explicitly highest weight vectors with respective reductive subalgebras.  
The two prominent examples
 in Introduction,
 i.e.
 the Rankin--Cohen bidifferential operators
 and the Juhl's conformally equivariant 
 differential operators,
 can be interpreted
 in the framework of the $F$-method
 as a special case of Theorem \ref{thm:1.4}.  
%See \cite{K08, K10} for some more multiplicity-free results.
\end{remark}

\section{A recipe of the $F$-method}

The idea of the $F$-method is to work on the branching problem of 
representations
 by taking
the Fourier transform of the nilpotent radical.
We shall explain this method in the complex setting
 where $H_{\mathbb{C}}$ is a parabolic subgroup $P_{\mathbb{C}}$
with abelian unipotent radical (see Theorem \ref{thm:surject} 
and Remark \ref{rem:cpx})
for simplicity.
A detaild proof will be given in \cite{KP}
 (see also \cite{KOSS}
 for a somewhat different formulation
 and normalization).

\subsection{Weyl algebra and algebraic Fourier transform}

Let $E$ be an $n$-dimensional vector space over $\C$. 
The Weyl algebra
$\mathcal D(E)$
is the ring of holomorphic differential operators on $E$
with polynomial coefficients. 

\begin{definition}[algebraic Fourier transform]\label{def:hat}
 We define an isomorphism of two Weyl algebras on $E$ and its dual space $E^\vee$:
\begin{equation}\label{eqn:algF}
\mathcal D(E)\to\mathcal D(E^\vee), \qquad T\mapsto \widehat T,
\end{equation}
which is induced by
\begin{equation}\label{eqn:Fgen}
\widehat{\frac\partial{\partial z_j}}:=-\zeta_j,\quad
\widehat z_j:=\frac\partial{\partial\zeta_j} \quad (1\leq j\leq n),
\end{equation} 
where $(z_1,\ldots,z_n)$ are coordinates on $E$ and $(\zeta_1,\ldots,\zeta_n)$ are 
the dual coordinates on $E^\vee$.
\end{definition}

\begin{remark}
(1)
The isomorphism \eqref{eqn:algF} is independent of the choice of coordinates.

(2)
An alternative way to get the isomorphism \eqref{eqn:algF}
or its variant is to use the Euclidean Fourier transform
$\mathcal{F}$ by choosing a real form $E(\mathbb{R})$ of $E$.
We then have
$$\widehat{\frac{\partial}{\partial z}}
= \sqrt{-1} \mathcal{F} \circ \frac{\partial}{\partial x} \circ \mathcal{F}^{-1},
\quad
\widehat{z} = - \sqrt{-1} \mathcal{F} \circ z \circ \mathcal{F}^{-1}
$$
as operators acting on the space $\mathcal{S}'(E^\vee)$
of Schwartz distributions.
This was the approach taken in \cite{KOSS}.
In particular,
 $\widehat T \ne {\mathcal{F}}\circ T \circ {\mathcal{F}}^{-1}$
 in our normalization here.  
The advantage of our normalization \eqref{eqn:Fgen} is 
 that  the commutative diagram in Theorem \ref{thm:313} 
 does not involve any power of $\sqrt{-1}$ that would otherwise depend on 
the degrees of differential operators.
As a consequence, the final step of the $F$-method  (see Step 5 below)
 as well as actual computations
 becomes simpler.
\end{remark}

\subsection{Infinitesimal action on principal series}

Let $\mathfrak{p} = \mathfrak{l} + \mathfrak{n}_+$
be a Levi decomposition of a parabolic subalgebra
 of $\mathfrak{g}$,
and 
$\mathfrak g=\mathfrak n_-+\mathfrak l+\mathfrak n_+$
the Gelfand--Naimark decomposition.
Since the following map
$$
\mathfrak n_-\times \mathfrak{l} \times \mathfrak n_+\to G_\C,\quad
(X,Z, Y)\mapsto (\exp X)(\exp Z)(\exp Y)
$$
is a local diffeomorphism near the origin,
we can define locally the projections
$p_-$ and $p_o$ 
from a neighbourhood of the identity
to the first and second factors $\mathfrak{n}_-$ and $\mathfrak{l}$,
respectively.
Consider the following two maps:
\begin{eqnarray*}
\alpha:\mathfrak g\times \mathfrak n_-\to\mathfrak l,&\quad&
 (Y,X)\mapsto \frac {d}{dt}\big\vert_{t=0}p_o\left(e^{tY}e^X\right),\\
\beta:\mathfrak g\times \mathfrak n_-\to\mathfrak n_-,&\quad&
 (Y,X)\mapsto \frac {d}{dt}\big\vert_{t=0}p_-\left(e^{tY}e^X\right).
\end{eqnarray*}
We may regard $\beta(Y,\cdot)$ as a vector field on $\mathfrak n_-$ by the
identification $\beta(Y,X)\in\mathfrak n_-\simeq T_X\mathfrak n_-.$

For $\mathfrak{l}$-module $\lambda$ on $V$,
we set $\mu:=\lambda^\vee\otimes\Lambda^{\dim}\mathfrak{n}_+$.
Since $\Lambda^{\dim\mathfrak{n}_+} \mathfrak{n}_+$
is one-dimensional,
we can and do identify the representation space with
$V^\vee$.
We inflate $\lambda$ and $\mu$ to $\mathfrak{p}$-modules
 by letting
$\mathfrak{n}_+$ act trivially.
Consider a Lie algebra homomorphism
\begin{equation}\label{eqn:dpi1}
d\pi_\mu:\mathfrak g\to\mathcal D(\mathfrak n_-)\otimes\operatorname{End}(V^\vee),
\end{equation}
defined for $F \in C^\infty(\mathfrak{n}_-,V^\vee)$ as
\begin{equation} \label{eqn:35}
\left(d\pi_\mu(Y)F\right)(X):=\mu(\alpha(Y,X))F(X)-(\beta(Y,\cdot\,)F)(X).
\end{equation}
If $(\mu,V^\vee)$ lifts to the parabolic subgroup
$P_{\mathbb{C}}$ of a reductive group $G_{\mathbb{C}}$
with Lie algebras $\mathfrak{p}$ and $\mathfrak{g}$ respectively,
then $d\pi_\mu$ is the differential representation of the induced
representation
$\operatorname{Ind}_{P_{\mathbb{C}}}^{G_{\mathbb{C}}}(V)$
(without $\rho$-shift).
We note
 that the Lie algebra homomorphism \eqref{eqn:35}
 is well-defined 
 without integrality condition of $\mu$.  
The $F$-method suggests to take
the algebraic Fourier transform \eqref{eqn:algF}
on the Weyl algebra $\mathcal D(\mathfrak n_-)$.
We then get another Lie algebra homomorphism
\begin{equation}\label{eqn:dpi2}
\widehat{d\pi_\mu}:\mathfrak g\to\mathcal D(\mathfrak n_+)\otimes\operatorname{End}(V^\vee).
\end{equation}
Then we have (see \cite{KP})
\begin{prop}
There is a natural isomorphism
\[
F_c: \operatorname{ind}_{\mathfrak{p}}^{\mathfrak{g}}(\lambda^\vee)
\overset{\sim}{\longrightarrow}
\operatorname{Pol}(\mathfrak{n}_+)\otimes V^\vee
\]
which intertwines the left\/ $\mathfrak{g}$-action on
$U(\mathfrak{g})\otimes_{U(\mathfrak{p})} V^\vee$
with $\widehat{d\pi_\mu}$.
\end{prop}

\subsection{Recipe of the $F$-method}\label{sec:44}

Our goal  is to find an explicit form of a $G'$-intertwining differential operator from $\mathcal V_X$ to $\mathcal W_Y$
in the upper right corner of Diagram 3.1 %~\ref{fig:1}.
Equivalently, what we call the $F$-method yields an explicit 
homomorphism belonging to
$\operatorname{Hom}_{\mathfrak g'}(\indpgprime(W^\vee),
\indpg(V^\vee)) \simeq
\operatorname{Hom}_{\mathfrak{p}'} (W^\vee,\operatorname{ind}_{\mathfrak{p}}^{\mathfrak{g}}(V^\vee))$
in the lower left corner of Diagram 3.1 %~\ref{fig:1}
 in the setting that 
 $\mathfrak{n}_+$ is abelian.

\medskip

The recipe of the $F$-method in this setting is stated
 as follows:
\begin{enumerate}
\item[Step 0.] Fix a finite dimensional representation $(\lambda,V)$ of $\mathfrak p=\mathfrak l+\mathfrak n_+$. 
\item[Step 1.] Consider a representation $\mu:=\lambda^\vee\,\otimes\, \Lambda^{\dim\mathfrak{n}_+}\mathfrak{n}_+$
 of the Lie algebra $\mathfrak p$.
Consider the restriction of the homomorphisms
 (\ref{eqn:dpi1}) and (\ref{eqn:dpi2}) 
to the subalgebra $\mathfrak{n}_+$:
\begin{eqnarray*}
d\pi_\mu:\mathfrak{n}_+&\to&\mathcal D(\mathfrak n_-)\otimes\operatorname{End}(V^\vee),\\
\widehat{d\pi_\mu}:\mathfrak{n}_+&\to&\mathcal D(\mathfrak n_+)\otimes\operatorname{End}(V^\vee).
\end{eqnarray*}

\item[Step 2.] Take a finite dimensional representation $W$ of the Lie algebra $\mathfrak p'$.
For the existence of nontrivial solutions in Step 3 below,
it is necessary and sufficient for $W$ to satisfy
\begin{equation}\label{eqn:Homnonzero}
\operatorname{Hom}_{\mathfrak{g}'}(\operatorname{ind}_{\mathfrak{p}'}^{\mathfrak{g}'}
(W^\vee),\indpg(V^\vee))\neq\{0\}.
\end{equation}
Choose $W$ satisfying \eqref{eqn:Homnonzero}
 if we know a priori an abstract
branching law of the restriction of
$\operatorname{ind}_{\mathfrak{p}}^{\mathfrak{g}}(V^\vee)$
to $\mathfrak{g}'$.
See \cite[Theorems 8.3 and 8.4]{K08} or \cite{K10} for some general formulae.  
Otherwise, we take $W$ to be any $\mathfrak{l}'$-irreducible component
of $S(\mathfrak{n}_+) \otimes V^\vee$
 and go to Step 3.  

\item[Step 3.]  Consider the system of 
partial differential equations for $\psi\in\operatorname{Pol}(\mathfrak n_+)\otimes 
V^\vee\otimes W$ which is $\mathfrak{l}'$-invariant
under the diagonal action:
 \begin{eqnarray}
\widehat{d\pi_\mu}( C) \psi=0 &&\mathrm{for}\, C\in\mathfrak n_+'.\label{eqn:Fmethod2}
\end{eqnarray}
Notice that  the equations
(\ref{eqn:Fmethod2}) are of second order.  
The solution space will be one-dimensional
 if we have chosen $W$ in Step 2
 such that
\begin{equation}
\label{eqn:Hom1}
\dim \operatorname{Hom}_{\mathfrak{g}'}(\operatorname{ind}_{\mathfrak{p}'}^{\mathfrak{g}'}
(W^\vee),\indpg(V^\vee))=1.  
\end{equation}

\item[Step 4.] Use invariant theory and 
reduce \eqref{eqn:Fmethod2} to another
system of differential equations on a lower dimensional space $S$.
Solve it. 

\item[Step 5.] Let $\psi$ be a polynomial solution to  (\ref{eqn:Fmethod2}) obtained in Step 4.
 Compute $(\operatorname {Symb}\otimes \operatorname{Id})^{-1}(\psi)$.
Here the symbol map
\[
\operatorname{Symb}:\operatorname{Diff}^{\operatorname{const}}({\mathfrak{n}}_-)
\overset \sim \to 
\operatorname{Pol}({\mathfrak{n}}_+)
\]
is a ring isomorphism
 given by the coordinates
\[
{\mathbb{C}}[\frac {\partial}{\partial z_1}, \cdots, \frac {\partial}{\partial z_n}]
\to 
{\mathbb{C}}[\xi_1, \cdots, \xi_n], 
\quad
\frac {\partial}{\partial z_j}
\mapsto 
\xi_j.
\]
\end{enumerate}
   
In case 
the Lie algebra representation $(\lambda,V)$ 
lifts to a group $P_\C$, we form a $G_\C$-equivariant holomorphic vector bundle $\mathcal V_{X_{\C}}$ over $X_\C=G_\C/P_\C$.
Likewise,
in case $W$ lifts to a group $P'_\C$,
 we form a $G'_\C$-equivariant holomorphic
vector bundle $\mathcal W_{Y_\C}$ over $Y_\C=G'_\C/P'_\C$.
Then
$(\operatorname{Symb} \otimes \operatorname{Id})^{-1}(\psi)$
in Step 5 gives an explicit formula of a $G'_{\mathbb{C}}$-equivariant differential
operator from $\mathcal{V}_{X_\C}$ to $\mathcal{W}_{Y_\C}$
in the coordinates of $\mathfrak n_-$
owing to Theorem \ref{thm:313} below.
This is what we wanted.  

\begin{remark}
In Step 2 we can find all such $W$ if we know a priori 
 (abstract) explicit branching laws.  
This is the case, 
 e.g., 
 in the setting of Theorem \ref{thm:1.4}.  
See Remark \ref{rem:Branching}.

 Conversely,
the differential equations in Step 3 sometimes give a useful information on
branching laws even when the restrictions are not completely reducible,
see \cite{KOSS}.
\end{remark}

For concrete constructions of equivariant differential operators
by using the $F$-method in various geometric settings,
we refer to \cite{KOSS,KP}.
A further application of the $F$-method to the construction of
\textit{non-local} operators will be discussed in another paper.
 
The key tool for the F-method is summarized as:
\begin{theorem}[\cite{KP}]\label{thm:313}
Let $P'_{\mathbb{C}}$ be a parabolic subgroup of $G'_{\mathbb{C}}$ compatible with 
a parabolic subgroup
$P_{\mathbb{C}}$ of $G_{\mathbb{C}}$. Assume further the nilradical\/ $\mathfrak n_+$ of\/ $\mathfrak{p}$ is abelian. 
Then the following diagram commutes:
%\begin{figure}[H]
%$$
\begin{align*}
&
\begin{matrix}
\operatorname{Hom}_\C (W^\vee,\indpg(V^\vee)) \,
\simeq
&\operatorname{Pol}(\mathfrak n_+)\otimes \operatorname{Hom}_\C(V,W)
& \stackrel{\operatorname{Symb}\otimes \mathrm{Id}}{\stackrel{\sim}{\longleftarrow}}\operatorname{Diff}^{\mathrm {const}}(\mathfrak n_-)\otimes \operatorname{Hom}_\C(V,W)\\
\cup & \circlearrowright& \cup\\
\operatorname{Hom}_{\mathfrak p'}(W^\vee,\indpg(V^\vee))
&\stackrel{D_{X_\C \to Y_\C}}{\longrightarrow}&
\operatorname{Diff}_{G_\C'}(\mathcal V_{X_\C},\mathcal W_{Y_\C}).
\end{matrix}
\\
&\hphantom{MMMMMMMMMMMMMMMM}\operatorname{Diagram 3.1} 
\end{align*}
%$$
%\caption{}
%\label{fig:1}
%\end{figure}
\end{theorem}

\end{document}